\documentclass[a4paper, 12pt]{article}
\usepackage{amssymb}

\usepackage{amsmath}
\usepackage{graphicx}

\input{tcilatex}

\begin{document}

\begin{center}
{\Huge Noncommutative recurrence over locally compact Hausdorff groups }

{\Huge \vspace{10mm} }{\normalsize Richard de Beer, Rocco Duvenhage\footnote{%
{\normalsize Corresponding author. \textit{E-mail address:}
rocco@postino.up.ac.za (R. Duvenhage).}} and Anton Str\"{o}h }

\emph{Department of Mathematics and Applied Mathematics}

\emph{University of Pretoria, 0002 Pretoria, South Africa}{\normalsize \ }

2005-4-8
\end{center}

\noindent\textbf{Abstract}

\bigskip

We extend previous results on noncommutative recurrence in unital $\ast $%
-algebras over the integers, to the case where one works over locally
compact Hausdorff groups. We derive a generalization of Khintchine's
recurrence theorem, as well as a form of multiple recurrence. This is done
using the mean ergodic theorem in Hilbert space, via the GNS construction.

\section{Introduction}

The simplest form of recurrence occurs in a dynamical system consisting of a
measure space $X$ with probability measure $\nu$, and a transformation $%
T:X\rightarrow X$ such that $\nu\left( T^{-1}(S)\right) =\nu(S)$ for all
measurable $S\subset X$. If $\nu(S)>0$ for some $S\subset X$, then there is
an $n\in\mathbb{N}=\{1,2,3,...\}$ such that $\nu\left( S\cap
T^{-n}(S)\right) >0$. This is essentially a pigeon hole principle for
measure spaces, and is usually referred to as Poincar\'{e}'s recurrence
theorem. Note that in this case the group over which we work is simply $%
\mathbb{Z}$. More precisely, we're working on the subset $\mathbb{N}$, since
we only consider $T^{-n}$ with $n\in\mathbb{N}$.

Recurrence theorems can also be studied in the noncommutative setting of
states on unital $\ast$-algebras and $C^{\ast}$-algebras, as done in [4] and
[2]. Typically the goal is to generalize existing measure theoretic
recurrence theorems to the $\ast$-algebraic setting. The measure theoretic
case is then recovered by taking the algebra to be a suitable commutative
algebra of measurable functions. The simplest way of doing this, is to use
the unital $\ast$-algebra $B_{\infty}(\Sigma)$ of bounded complex-valued
measurable functions on the measure space $X$, with $\Sigma$ the $\sigma$%
-algebra of measurable sets of $X$. The state $\omega$ on the algebra is
simply $\omega(f):=\int_{X}fd\nu$, so the state represents the measure,
while $T$ is represented on the algebra by the Koopman construction $%
\tau(f)=f\circ T$. Note that $\omega(\chi_{S})=\nu(S)$ and $%
\tau(\chi_{S})=\chi_{T^{-1}(S)}$. By following this type of recipe, the
results in Sections 3 and 4 can be applied to the measure theoretic case.

In this paper we continue our work in [2]. We study recurrence in unital $%
\ast$-algebras as before, but instead of just working over the group $%
\mathbb{Z}$ as mentioned above, we will consider locally compact Hausdorff
groups and suitable subsets thereof, namely subsemigroups (similar to the
subset $\mathbb{N}$ of $\mathbb{Z}$). The main result is an extension of
Khintchine's recurrence theorem in Section 3, which is subsequently used to
prove a simple multiple recurrence result in Section 4. The latter result is
inspired by the work of Furstenberg [3] on extensions of Poincar\'{e}
recurrence to recurrence theorems of the form 
\begin{equation*}
\nu\left( S\cap T^{n}S\cap T^{2n}S\cap...\cap T^{(k-1)n}S\right) >0 
\end{equation*}
which Furstenberg used to give an alternative proof of Szemer\'{e}di's
theorem in combinatorial number theory. Our result is of a different form
than Furstenberg's, however.

The main tool we use, is the mean ergodic theorem in Hilbert space, which we
review in Section 2. A Hilbert space version of Khintchine's recurrence
theorem is also proved in Section 2, and applied in the subsequent sections.
In the recurrence theorems we need to make stronger assumptions than in the
mean ergodic theorem, namely that the subsemigroup of the group over which
we work is abelian, or that the group itself is unimodular.

\section{Recurrence in a Hilbert space setting}

We start with a review of the mean ergodic theorem. This is based on
Petersen [5] and also Bratteli and Robinson [1]. The former discusses the
theorem over the group $G=\mathbb{Z}$ (see the more general form below),
while the latter gives it in an abstract form involving no group.

First, consider a function $f:G\rightarrow\mathfrak{H}$ where $G$ is a
locally compact Hausdorff group with right Haar measure $\mu$, and $%
\mathfrak{H}$ a complex Hilbert space, such that $G\ni g\mapsto\left\langle
f(g),x\right\rangle $ is Borel measurable for every $x\in\mathfrak{H}$.
We'll take the second slot in the inner product to be the linear one. If $%
\Lambda\subset G$ is Borel with $\mu(\Lambda)<\infty$, and $f$ is bounded on 
$\Lambda$, say $\left| \left| f(g)\right| \right| \leq b$ for all $%
g\in\Lambda$ for some positive $b\in\mathbb{R}$, then we can define $%
\int_{\Lambda}fd\mu$ by $\left\langle \int_{\Lambda}fd\mu,x\right\rangle
:=\int_{\Lambda}\left\langle f(g),x\right\rangle d\mu(g)$ for all $x\in%
\mathfrak{H}$ using the Riesz representation theorem. So we also have $%
\left\langle x,\int_{\Lambda}fd\mu\right\rangle =\int_{\Lambda}\left\langle
x,f(g)\right\rangle d\mu(g)$. We will also use the notation $\int_{\Lambda
}f(g)dg=\int_{\Lambda}fd\mu$. One can then easily prove all the standard
properties for this integral, like linearity and%
\begin{equation}
\int_{\Lambda}f(gh)dg=\int_{\Lambda h}fd\mu   \tag{1}
\end{equation}%
\begin{equation}
\int_{\Lambda}Af(g)dg=A\int_{\Lambda}fd\mu   \tag{2}
\end{equation}%
\begin{equation}
\int_{\Lambda}xd\mu=x\int_{\Lambda}d\mu=\mu(\Lambda)x   \tag{3}
\end{equation}%
\begin{equation}
\int_{\Lambda_{1}\cup\Lambda_{2}}fd\mu=\int_{\Lambda_{1}}fd\mu+\int
_{\Lambda_{2}}fd\mu   \tag{4}
\end{equation}
and%
\begin{equation}
\left\| \int_{\Lambda}fd\mu\right\| \leq b\mu(\Lambda)   \tag{5}
\end{equation}
for every $h\in G$, $A\in B(\mathfrak{H})$, $x\in\mathfrak{H}$ and Borel $%
\Lambda_{1},\Lambda_{2}\subset G$ of finite measure on which $f$ is bounded,
with $\mu(\Lambda_{1}\cap\Lambda_{2})=0$, where $B(\mathfrak{H})$ denotes
the algebra of bounded linear operators $\mathfrak{H}\rightarrow\mathfrak{H}$%
. We will use these properties in the sequel.

A \textit{net} is family $\left\{ \Lambda_{\alpha}\right\} $ of subsets of $G
$ indexed by a directed set. If a $K\subset G$ (with equality allowed) has
the property that $gh\in K$ for all $g,h\in K$, we'll call $K$ a \textit{%
subsemigroup} of $G$. We call a net $\left\{ \Lambda_{\alpha }\right\} $ of
Borel subsets of $G$ \textit{space-filling in} $K$ if $\Lambda_{\alpha}%
\subset K$, $\mu(\Lambda_{\alpha})<\infty$, and 
\begin{equation*}
\lim_{\alpha}\frac{\mu\left( \Lambda_{\alpha}\Delta\left( \Lambda_{\alpha
}g\right) \right) }{\mu\left( \Lambda_{\alpha}\right) }=0 
\end{equation*}
for all $g\in K$, where we assume that $\mu(\Lambda_{\alpha})>0$ for $\alpha$
large enough (i.e. $\alpha\geq\alpha_{0}$ for some $\alpha_{0}$). Here $%
A\Delta B:=\left( A\cup B\right) \backslash\left( A\cap B\right) $. Now we
can state

\bigskip

\noindent\textbf{The mean ergodic theorem. }\textit{Let $G$ be a locally
compact Hausdorff group with right Haar measure $\mu$, and consider a Borel
measurable subsemigroup $K$ of $G$. Let $U:K\rightarrow B(\mathfrak{H}%
):g\mapsto U_{g}$ be such that $\left\| U_{g}\right\| \leq1$, $%
U_{g}U_{h}=U_{gh}$ for all $g,h\in K$, and $K\ni g\mapsto\left\langle
U_{g}x,y\right\rangle $ is Borel measurable for all $x,y\in\mathfrak{H}$.
Take $P$ to be the projection of $\mathfrak{H}$ onto $V:=\left\{ x\in 
\mathfrak{H}:U_{g}x=x\text{ for all }g\in K\right\} $. For any space-filling
net $\left\{ \Lambda_{\alpha}\right\} $ in $K$ we then have 
\begin{equation*}
\lim_{\alpha}\frac{1}{\mu\left( \Lambda_{\alpha}\right) }\int_{\Lambda
_{\alpha}}U_{g}xdg=Px 
\end{equation*}
for all $x\in\mathfrak{H}$.}

\bigskip

\noindent\textbf{Proof. }Set $N:=\overline{\text{span}\left\{ x-U_{g}x:x\in%
\mathfrak{H},g\in K\right\} }$. For any $g$, a fixed point of $U_{g}^{\ast}$
is a fixed point of $U_{g}$, and vice versa, since $\left| \left|
U_{g}^{\ast}\right| \right| \leq1$. From this it follows that $V=N^{\bot}$,
which means in particular that $V$ is a closed subspace of $\mathfrak{H}$.
Set 
\begin{equation*}
I_{\alpha}(x):=\frac{1}{\mu(\Lambda_{\alpha})}\int_{\Lambda_{\alpha}}U_{g}xdg%
\text{ \ \ .}
\end{equation*}
We now prove that $\lim_{\alpha}I_{\alpha}(x)=0$ for $x\in N$. First let $%
x=y-U_{h}y$ for some $y\in\mathfrak{H}$ and $h\in K$, then we have from (1)
and (4) that 
\begin{align*}
I_{\alpha}(x) & =\frac{1}{\mu(\Lambda_{\alpha})}\int_{\Lambda_{%
\alpha}}U_{g}ydg-\frac{1}{\mu(\Lambda_{\alpha})}\int_{\Lambda_{%
\alpha}h}U_{g}ydg \\
& =\frac{1}{\mu(\Lambda_{\alpha})}\int_{\Lambda_{\alpha}\backslash\left(
\Lambda_{\alpha}\cap(\Lambda_{\alpha}h)\right) }U_{g}ydg-\frac{1}{\mu
(\Lambda_{\alpha})}\int_{(\Lambda_{\alpha}h)\backslash\left( \Lambda_{\alpha
}\cap(\Lambda_{\alpha}h)\right) }U_{g}ydg
\end{align*}
hence 
\begin{align*}
\left| \left| I_{\alpha}(x)\right| \right| & \leq\frac{1}{\mu
(\Lambda_{\alpha})}\left| \left| \int_{\Lambda_{\alpha}\backslash\left(
\Lambda_{\alpha}\cap(\Lambda_{\alpha}h)\right) }U_{g}ydg\right| \right| +%
\frac{1}{\mu(\Lambda_{\alpha})}\left| \left| \int_{(\Lambda_{\alpha
}h)\backslash\left( \Lambda_{\alpha}\cap(\Lambda_{\alpha}h)\right)
}U_{g}ydg\right| \right| \\
& \leq\left| \left| y\right| \right| \frac{\mu\left( \Lambda_{\alpha
}\Delta(\Lambda_{\alpha}h)\right) }{\mu(\Lambda_{\alpha})}
\end{align*}
by (5), since $\left| \left| U_{g}\right| \right| \leq1$, so $\lim
_{\alpha}I_{\alpha}(x)=0$. However, we need this for any $x\in N$, so set $%
N_{0}:=\left\{ y-U_{g}y:y\in\mathfrak{H},g\in K\right\} $. Then for any $%
\varepsilon>0$ there is a $y\in$ span$N_{0}$ such that $\left| \left|
x-y\right| \right| <\varepsilon$, say $y=\sum_{j=1}^{m}x_{j}$ where $%
x_{j}\in N_{0}$. Therefore%
\begin{equation*}
\left| \left\| I_{\alpha}(x)\right\| -\left\| I_{\alpha}(y)\right\| \right|
\leq\left| \left| I_{\alpha}(x)-I_{\alpha}(y)\right| \right| \leq\frac{1}{%
\mu(\Lambda_{\alpha})}\left| \left| x-y\right| \right|
\int_{\Lambda_{\alpha}}d\mu<\varepsilon 
\end{equation*}
while 
\begin{equation*}
\left| \left| I_{\alpha}(y)\right| \right| \leq\sum_{j=1}^{m}\left| \left|
I_{\alpha}(x_{j})\right| \right| \rightarrow0 
\end{equation*}
in the $\alpha$ limit, as shown above. Hence $\lim_{\alpha}I_{\alpha}(x)=0$.

For any $x\in\mathfrak{H}$, write $x=x_{0}+Px$, where $x_{0}=(1-P)x\in
V^{\bot}=N$, then 
\begin{equation*}
\left| \left| I_{\alpha}(x)-Px\right| \right| =\left| \left| I_{\alpha
}(x_{0})+I_{\alpha}(Px)-Px\right| \right| =\left| \left|
I_{\alpha}(x_{0})\right| \right| \rightarrow0 
\end{equation*}
in the $\alpha$ limit, since $I_{a}(Px)=\frac{1}{\mu(\Lambda_{\alpha})}%
\int_{\Lambda_{\alpha}}Pxd\mu=Px$ by the definition of $P$ and (3). $\square$

\bigskip

Using this theorem we can prove a Hilbert space version of Khintchine's
recurrence theorem:

\bigskip

\noindent\textbf{Theorem 2.1.}{\textit{\ Consider the situation given in the
mean ergodic theorem above, but assume $K$ is abelian, i.e. $gh=hg$ for all $%
g,h\in K$. Take any $x,y\in\mathfrak{H}$ and $\varepsilon>0$. Then there is
an $\alpha_{0}$ such that%
\begin{equation*}
\left| \frac{1}{\mu(\Lambda_{\alpha_{0}})}\int_{\Lambda_{\alpha_{0}}h}\left%
\langle x,U_{g}y\right\rangle dg\right| >\left| \left\langle
x,Py\right\rangle \right| -\varepsilon 
\end{equation*}
for all $h\in K$. In particular, for every $h\in K$ there is a $g\in
\Lambda_{\alpha_{0}}h$ such that%
\begin{equation*}
\left| \left\langle x,U_{g}y\right\rangle \right| >\left| \left\langle
x,Py\right\rangle \right| -\varepsilon\text{ \ \ .}
\end{equation*}
}}

\bigskip

\noindent\textbf{Proof. }By the mean ergodic theorem there is an $\alpha_{0}$
such that 
\begin{equation*}
\left| \left| \frac{1}{\mu(\Lambda_{\alpha_{0}})}\int_{\Lambda_{%
\alpha_{0}}}U_{g}ydg-Py\right| \right| <\frac{\varepsilon}{\left| \left|
x\right| \right| +1}
\end{equation*}
while by definition of $P$ we have $U_{h}Py=Py$ for all $h\in K$. Using
these two facts along with (1) and (2) and the fact that $K$ is abelian, we
get 
\begin{align*}
\left| \left| \frac{1}{\mu(\Lambda_{\alpha_{0}})}\int_{\Lambda_{%
\alpha_{0}}h}U_{g}ydg-Py\right| \right| & =\left| \left| \frac{1}{\mu
(\Lambda_{\alpha_{0}})}\int_{\Lambda_{\alpha_{0}}}U_{gh}ydg-Py\right| \right|
\\
& =\left| \left| \frac{1}{\mu(\Lambda_{\alpha_{0}})}U_{h}\int
_{\Lambda_{\alpha_{0}}}U_{g}ydg-U_{h}Py\right| \right| \\
& \leq\frac{\varepsilon}{\left| \left| x\right| \right| +1}
\end{align*}
since $\left| \left| U_{h}\right| \right| \leq1$. Hence%
\begin{equation*}
\left| \frac{1}{\mu(\Lambda_{\alpha_{0}})}\int_{\Lambda_{\alpha_{0}}h}\left%
\langle x,U_{g}y\right\rangle dg-\left\langle x,Py\right\rangle \right|
=\left| \left\langle x,\frac{1}{\mu(\Lambda_{\alpha_{0}})}%
\int_{\Lambda_{\alpha_{0}}h}U_{g}ydg-Py\right\rangle \right| <\varepsilon 
\end{equation*}
from which the result follows. $\square$

\bigskip

In this result we assumed $K$ to be abelian, but if we assume $G$ is
unimodular, i.e. its right Haar measure is also a left Haar measure, then
essentially the same proof also works for non-abelian $K$ to give

\bigskip

\noindent\textbf{Theorem 2.2. }\textit{Consider the situation given in the
mean ergodic theorem above, but assume $G$ is unimodular. For any $x,y\in%
\mathfrak{H}$ and $\varepsilon>0$ there then exists an $\alpha_{0}$ such
that 
\begin{equation*}
\left| \frac{1}{\mu(\Lambda_{\alpha_{0}})}\int_{h\Lambda_{\alpha_{0}}}\left%
\langle x,U_{g}y\right\rangle dg\right| >\left| \left\langle
x,Py\right\rangle \right| -\varepsilon 
\end{equation*}
for all $h\in K$. $\square$}

\bigskip

This works, since if $\mu$ is also a left Haar measure, one has $\int
_{\Lambda}f(hg)dg=\int_{h\Lambda}fd\mu$ similar to (1).

\section{$\ast$-dynamical systems}

Let $L(V)$ denote the space of all linear operators $V\rightarrow V$ with $V$
a vector space.

\bigskip

\noindent\textbf{Definition 3.1. }Let $\omega$ be a state on a unital $\ast $%
-algebra $\mathfrak{A}$, let $G$ be a locally compact Hausdorff group with
right Haar measure $\mu$, and $K$ a Borel measurable subsemigroup of $G$.
Consider a $\tau:K\rightarrow L(\mathfrak{A}):g\mapsto\tau_{g}$ with 
\begin{align*}
\tau_{g}\circ\tau_{h} & =\tau_{gh} \\
\tau_{g}(1) & =1 \\
\omega\left( \tau_{g}(A)^{\ast}\tau_{g}(A)\right) & \leq\omega(A^{\ast}A)
\end{align*}
for all $g,h\in K$ and $A\in\mathfrak{A}$, and $K\ni g\mapsto\omega\left(
A^{\ast}\tau_{g}(B)\right) $ Borel measurable for all $A,B\in\mathfrak{A}$.
Then we'll call $\left( \mathfrak{A},\omega,\tau,K\right) $ a $\ast $\textit{%
-dynamical system}. $\square$

\bigskip

Given a state $\omega$ on a unital $\ast$-algebra $\mathfrak{A}$, the GNS
construction provides us with a cyclic representation $\left( \mathfrak{G}%
,\pi,\Omega\right) $ where $\mathfrak{G}$ is an inner product space, $\pi:%
\mathfrak{A}\rightarrow L(\mathfrak{G})$ is linear, and $\iota :\mathfrak{A}%
\rightarrow\mathfrak{G}:A\mapsto\pi(A)\Omega$ is surjective. Also, $%
\omega(A^{\ast}B)=\left\langle \iota(A),\iota(B)\right\rangle $ for all $%
A,B\in\mathfrak{A}$. Then for $\tau_{g}$ as above 
\begin{equation*}
U_{g}:\mathfrak{G}\rightarrow\mathfrak{G}:\iota(A)\mapsto\iota(\tau_{g}(A)) 
\end{equation*}
is well-defined, linear, and $\left| \left| U_{g}\right| \right| \leq1$. We
can therefore uniquely extend $U_{g}$ to the completion $\mathfrak{H}$ of $%
\mathfrak{G}$, such that $\left| \left| U_{g}\right| \right| \leq1$. Call $%
g\mapsto U_{g}$ the GNS representation of $\tau$.

\bigskip

\noindent\textbf{Proposition 3.2. }\textit{For a $\ast$-dynamical system $%
\left( \mathfrak{A},\omega,\tau,K\right) $, the GNS representation $U$ of $%
\tau$ on a Hilbert space $\mathfrak{H}$ has the following properties: $%
U_{g}U_{h}=U_{gh}$ for all $g,h\in K$, and $K\in g\mapsto\left\langle
x,U_{g}y\right\rangle $ is Borel measurable for all $x,y\in\mathfrak{H}$.}

\bigskip

\noindent\textbf{Proof.} For $A\in\mathfrak{A}$ one has $U_{g}U_{h}%
\iota(A)=U_{g}\iota\left( \tau_{h}(A)\right) =\iota\left( \tau_{g}\left(
\tau_{h}(A)\right) \right) =\iota\left( \tau_{gh}(A)\right) =U_{gh}\iota(A)$%
, and by continuity of $U_{g}$ on $\mathfrak{H}$, this extends to $%
U_{g}U_{h}x=U_{gh}x$ for all $x\in\mathfrak{H}$. By the definition of a $\ast
$-dynamical system, $g\mapsto\omega\left( A^{\ast}\tau_{g}(B)\right)
=\left\langle \iota(A),U_{g}\iota(B)\right\rangle $ is Borel measurable, and
since the pointwise limit of a sequence of measurable functions is
measurable, we need only consider $\left\langle
x_{n},U_{g}y_{n}\right\rangle \rightarrow\left\langle x,U_{g}y\right\rangle $
where $x,y\in\mathfrak{H}$ and $x_{n},y_{n}\in\mathfrak{G}$ (as defined
above) such that $x_{n}\rightarrow x$ and $y_{n}\rightarrow y$, keeping in
mind that $x_{n}=\iota(A_{n})$ and $y_{n}=\iota(B_{n})$ for some $%
A_{n},B_{n}\in\mathfrak{A}$. $\square$

\bigskip

Now we can state a recurrence theorem for $\ast$-dynamical systems,
containing in particular the conventional form of the Khintchine recurrence
theorem (which includes the measure theoretic version over $K=\mathbb{N}$,
as a special case; see Petersen [5]):

\bigskip

\noindent\textbf{Theorem 3.3. }\textit{Let $\left( \mathfrak{A},\omega
,\tau,K\right) $ be a $\ast$-dynamical system, but assume $K$ is abelian.
Let $\left\{ \Lambda_{\alpha}\right\} $ be a space-filling net in $K$. Then
for any $A,B\in\mathfrak{A}$ and $\varepsilon>0$, there exists an $\alpha_{0}
$ such that%
\begin{equation*}
\left| \frac{1}{\mu(\Lambda_{\alpha_{0}})}\int_{\Lambda_{\alpha_{0}}h}\omega%
\left( A^{\ast}\tau_{g}(B)\right) dg\right| >\left| \lim_{\alpha }\frac{1}{%
\mu(\Lambda_{\alpha})}\int_{\Lambda_{\alpha}}\omega\left( A^{\ast
}\tau_{g}(B)\right) dg\right| -\varepsilon 
\end{equation*}
for all $h\in K$. In particular, if $A=B$, then for every $h\in K$ 
\begin{equation*}
\left| \omega\left( A^{\ast}\tau_{g}(A)\right) \right| >\left|
\omega(A)\right| ^{2}-\varepsilon 
\end{equation*}
for some $g\in\Lambda_{\alpha_{0}}h$.}

\bigskip

\noindent\textbf{Proof.} We use the GNS\ construction discussed above to
represent $\tau$ by $U$, and set $x:=\iota(A)$ and $y:=\iota(B)$. By the
mean ergodic theorem 
\begin{equation*}
\left\langle x,Py\right\rangle =\lim_{\alpha}\frac{1}{\mu(\Lambda_{\alpha})}%
\int_{\Lambda_{\alpha}}\left\langle x,U_{g}y\right\rangle dg=\lim_{\alpha }%
\frac{1}{\mu(\Lambda_{\alpha})}\int_{\Lambda_{\alpha}}\omega\left( A^{\ast
}\tau_{g}(B)\right) dg\text{ \ \ .}
\end{equation*}
The first part of the result now follows immediately from Theorem 2.1. The
second part also follows, since $\left| \omega(A)\right| =\left|
\omega(1^{\ast}A)\right| =\left| \left\langle \iota(1),\iota(A)\right\rangle
\right| =\left| \left\langle \Omega,x\right\rangle \right| =\left|
\left\langle P\Omega,x\right\rangle \right| =\left| \left\langle
\Omega,Px\right\rangle \right| \leq\left| \left| \Omega\right| \right|
\left| \left| Px\right| \right| =\left| \left| Px\right| \right| =\sqrt{%
\left\langle x,Px\right\rangle }$, where we've used the fact that $%
P\Omega=\Omega$, which follows from $U_{g}\Omega=U_{g}\iota(1)=\iota(\tau
_{g}(1))=\iota(1)=\Omega$. $\square$

\bigskip

Even though $K$ has to a be abelian in this theorem, one could still have a $%
\ast$- dynamical system $\left( \mathfrak{A},\omega,\tau,H\right) $ with $H$
non-abelian, and then apply the theorem to various abelian $K\subset H$ with 
$K$ a Borel measurable subsemigroup of the underlying group $G$.

\bigskip

\noindent\textbf{Remarks on ergodicity. }Call a $\ast$-dynamical system $%
\left( \mathfrak{A},\omega,\tau,K\right) $ \textit{ergodic} when 
\begin{equation*}
\lim_{\alpha}\frac{1}{\mu(\Lambda_{\alpha})}\int_{\Lambda_{\alpha}}\omega%
\left( A\tau_{g}(B)\right) dg=\omega(A)\omega(B) 
\end{equation*}
for all $A,B\in\mathfrak{A}$ and some space-filling net $\left\{
\Lambda_{\alpha}\right\} $ in $K$.

In the GNS representation, and with $P$ as in the mean ergodic theorem, the
above definition of ergodicity is equivalent to $P$ having a one dimensional
range, and in particular the definition is independent of which
space-filling net in $K$ is used. In fact, $P=\Omega\otimes\Omega$ in case
of ergodicity, where ($x\otimes y)z:=x\left\langle y,z\right\rangle $ for
all $x,y,z\in \mathfrak{H}$. We see this as follows:

Since $P\Omega=\Omega$ as we saw in the proof of Theorem 3.3, it follows
that $P$ having one-dimensional range is equivalent to $P=\Omega\otimes\Omega
$. Now, if $P=\Omega\otimes\Omega$, then the mean ergodic theorem tells us
that 
\begin{align*}
\lim_{\alpha}\frac{1}{\mu(\Lambda_{\alpha})}\int_{\Lambda_{\alpha}}\omega%
\left( A\tau_{g}(B)\right) dg & =\lim_{\alpha}\frac{1}{\mu (\Lambda_{\alpha})%
}\int_{\Lambda_{\alpha}}\left\langle
\iota(A^{\ast}),U_{g}\iota(B)\right\rangle dg \\
& =\left\langle \iota(A^{\ast}),P\iota(B)\right\rangle \\
& =\left\langle \iota(A^{\ast}),(\Omega\otimes\Omega)\iota(B)\right\rangle \\
& =\left\langle \iota(A^{\ast}),\Omega\right\rangle \left\langle \Omega
,\iota(B)\right\rangle \\
& =\omega(A)\omega(B)\text{ \ \ .}
\end{align*}
Conversely, if the system is ergodic, a similar argument shows that $%
\left\langle \iota(A^{\ast}),P\iota(B)\right\rangle =\omega(A)\omega
(B)=\left\langle \iota(A^{\ast}),\Omega\right\rangle \left\langle \Omega
,\iota(B)\right\rangle =\left\langle \iota(A^{\ast}),(\Omega\otimes
\Omega)\iota(B)\right\rangle $, and since $\mathfrak{G}$ is dense in $%
\mathfrak{H}$, it follows that $P=\Omega\otimes\Omega$.

A corollary of Theorem 3.3 for ergodic systems is clearly

\bigskip

\noindent\textbf{Corollary 3.4. }\textit{Assume $\left( \mathfrak{A}%
,\omega,\tau,K\right) $ given in Theorem 3.3 is ergodic. For any $A,B\in%
\mathfrak{A}$ and $\varepsilon>0$, there then exists an $\alpha_{0}$ such
that 
\begin{equation*}
\left| \frac{1}{\mu(\Lambda_{\alpha_{0}})}\int_{\Lambda_{\alpha_{0}}h}\omega%
\left( A\tau_{g}(B)\right) dg\right| >\omega(A)\omega(B)-\varepsilon 
\end{equation*}
for all $h\in K$. $\square$}

\bigskip

While ergodicity can be formulated, and proven equivalent to $P$ having
one-dimensional range, even when $K$ is not abelian, as we did above,
Theorem 3.3 and Corollary 3.4 are only stated for abelian $K$, though $G$ is
allowed to be non-abelian. However, using Theorem 2.2, Theorem 3.3 can be
modified to

\bigskip

\noindent\textbf{Theorem 3.5. }\textit{Let $\left( \mathfrak{A},\omega
,\tau,K\right) $ be a $\ast$-dynamical system, but assume\ that the
underlying group $G$ is unimodular. Let $\left\{ \Lambda_{\alpha}\right\} $
be a space-filling net in $K$. Then for any $A,B\in\mathfrak{A}$ and $%
\varepsilon>0$, there exists an $\alpha_{0}$ such that%
\begin{equation*}
\left| \frac{1}{\mu(\Lambda_{\alpha_{0}})}\int_{h\Lambda_{\alpha_{0}}}\omega%
\left( A^{\ast}\tau_{g}(B)\right) dg\right| >\left| \lim_{\alpha }\frac{1}{%
\mu(\Lambda_{\alpha})}\int_{\Lambda_{\alpha}}\omega\left( A^{\ast
}\tau_{g}(B)\right) dg\right| -\varepsilon 
\end{equation*}
for all $h\in K$. $\square$}

\bigskip

We can modify Corollary 3.4 in a corresponding way.

\section{Towards multiple recurrence}

In this section we study a form of multiple recurrence, inspired by
Furstenberg's work, as mentioned in the introduction. Also refer to Petersen
[5] for a discussion of multiple recurrence in the measure theoretic setting
over the group $G=\mathbb{Z}$. We will formulate our results for an abelian
subsemigroup $K$ of a locally compact Hausdorff group $G$, but as with
Theorem 3.5, the results in this section can easily be modified to the case
where $G$ is unimodular, but $K$ not necessarily abelian.

Let $\mathfrak{A}\otimes\mathfrak{B}$ denote the algebraic tensor product of
the $\ast$-algebras $\mathfrak{A}$ and $\mathfrak{B}$. First we state a
result about simultaneous recurrence in more than one system:

\bigskip

\noindent\textbf{Proposition 4.1. }\textit{Let $\left( \mathfrak{A}%
_{j},\omega_{j},\tau_{j},K\right) $ be a $\ast$-dynamical system such that $%
\omega\left( \tau_{j,g}(A)^{\ast}\tau_{j,g}(B)\right) =\omega(A^{\ast}B)$
for all $A,B\in\mathfrak{A}_{j}$ and $g\in K$, for $j=1,...,q$, and assume $K
$ is abelian. (Here we use the notation $\tau_{j}:K\rightarrow L(\mathfrak{A}%
):g\mapsto\tau_{j,g}$.) Let $\left\{ \Lambda_{\alpha}\right\} $ be a
space-filling net in $K$. Then for any $A_{j},B_{j}\in\mathfrak{A}_{j}$ and $%
\varepsilon>0$, there exists an $\alpha_{0}$ such that%
\begin{align*}
& \left| \frac{1}{\mu(\Lambda_{\alpha_{0}})}\int_{\Lambda_{\alpha_{0}}h}%
\omega_{1}\left( A_{1}^{\ast}\tau_{1,g}(B_{1})\right) ...\omega_{q}\left(
A_{q}^{\ast}\tau_{q,g}(B_{q})\right) dg\right| \\
& >\left| \lim_{\alpha}\frac{1}{\mu(\Lambda_{\alpha})}\int_{\Lambda_{\alpha
}}\omega_{1}\left( A_{1}^{\ast}\tau_{1,g}(B_{1})\right) ...\omega_{q}\left(
A_{q}^{\ast}\tau_{q,g}(B_{q})\right) dg\right| -\varepsilon
\end{align*}
for all $h\in K$. In particular, if $A_{j}=B_{j}$, then for every $h\in K$%
\begin{equation*}
\left| \omega_{1}\left( A_{1}^{\ast}\tau_{1,g}(A_{1})\right) ...\omega
_{q}\left( A_{q}^{\ast}\tau_{q,g}(A_{q})\right) \right| >\left| \omega
_{1}\left( A_{1}\right) ...\omega_{q}\left( A_{q}\right) \right|
^{2}-\varepsilon 
\end{equation*}
for some $g\in\Lambda_{\alpha_{0}}h$.}

\bigskip

\noindent\textbf{Proof.} Set $\mathfrak{A}:=\mathfrak{A}_{1}\otimes
...\otimes\mathfrak{A}_{q}$, $\omega:=\omega_{1}\otimes...\otimes\omega_{q}$
and $\tau_{g}:=\tau_{1,g}\otimes...\otimes\tau_{q,g}$. We first show that $%
\left( \mathfrak{A},\omega,\tau,K\right) $ is a $\ast$-dynamical system. (It
is in this step that we require the additional condition $\omega\left(
\tau_{j,g}(A)^{\ast}\tau_{j,g}(B)\right) =\omega(A^{\ast}B)$.) First
consider $A=A_{1}\otimes A_{2}$ and $B=B_{1}\otimes B_{2}$ where $%
A_{j},B_{j}\in\mathfrak{A}_{j}$. Then $\left(
\tau_{1,g}\otimes\tau_{2,g}\right) \left( \left(
\tau_{1,h}\otimes\tau_{2,h}\right) (A)\right) =\left(
\tau_{1,g}\otimes\tau_{2,g}\right) \left( \tau_{1,h}(A_{1})\otimes\tau
_{2,h}(A_{2})\right) =\left( \tau_{1,gh}\otimes\tau_{2,gh}\right) (A)$, $%
\left( \tau_{1,g}\otimes\tau_{2,g}\right) (1)=\left( \tau_{1,g}\otimes
\tau_{2,g}\right) (1\otimes1)=1\otimes1=1$, and%
\begin{align*}
& \omega_{1}\otimes\omega_{2}\left( \left[ \left( \tau_{1,g}\otimes
\tau_{2,g}\right) (A)\right] ^{\ast}\left( \tau_{1,g}\otimes\tau
_{2,g}\right) (B)\right) \\
& =\omega_{1}\left( \tau_{1,g}(A_{1})^{\ast}\tau_{1,g}(B_{1})\right)
\omega_{2}\left( \tau_{2,g}(A_{2})^{\ast}\tau_{2,g}(B_{2})\right) \\
& =\omega_{1}(A_{1}^{\ast}B_{1})\omega_{2}(A_{2}^{\ast}B_{2}) \\
& =\omega_{1}\otimes\omega_{2}(A^{\ast}B)
\end{align*}
Furthermore, 
\begin{equation*}
K\ni g\mapsto\omega_{1}\otimes\omega_{2}\left( A^{\ast}\left( \tau
_{1,g}\otimes\tau_{2,g}\right) (B)\right) =\omega_{1}\left( A_{1}{}^{\ast
}\tau_{1,g}(B_{1})\right) \omega_{2}\left(
A_{2}{}^{\ast}\tau_{2,g}(B_{2})\right) 
\end{equation*}
is Borel measurable. All these facts then also hold for any $A,B\in 
\mathfrak{A}_{1}\otimes\mathfrak{A}_{2}$, since these elements have the form 
$A=\sum_{k=1}^{M}A_{1,k}\otimes A_{2,k}$ where $A_{j,k}\in\mathfrak{A}_{j}$.
By induction this can be extended to obtain $\tau_{g}\circ\tau_{h}=\tau_{gh}$%
, $\tau_{g}(1)=1$, $\omega\left( A^{\ast}\tau_{g}(B)\right) =\omega(A^{\ast
}B)$, and $K\ni g\mapsto\omega\left( A^{\ast}\tau_{g}(B)\right) $ Borel
measurable. In particular, $\left( \mathfrak{A},\omega,\tau,K\right) $ is a $%
\ast$-dynamical system. Applying Theorem 3.3 to $A:=A_{1}\otimes...\otimes
A_{q}$ and $B:=B_{1}\otimes...\otimes B_{q}$, the proposition is proved. $%
\square$

\bigskip

We are now going to apply this result to prove a form of multiple
recurrence. Given a $\tau:K\rightarrow L(\mathfrak{A})$ such that $%
\tau_{g}\circ\tau _{h}=\tau_{gh}$, and we want to construct a $%
\sigma:K\rightarrow L(\mathfrak{A}):g\mapsto\tau_{\varphi(g)}$ such that we
also have $\sigma _{g}\circ\sigma_{h}=\sigma_{gh}$, then it seems sensible
to take $\varphi:K\rightarrow K$ such that $\varphi(g)\varphi(h)=\varphi(gh)$%
. Such $\varphi$'s will determine the pattern of the multiple recurrence
(see Corollary 4.3 and the discussion following it).

\bigskip

\noindent\textbf{Theorem 4.2. }\textit{Consider a $\ast$-dynamical system $%
\left( \mathfrak{A},\omega,\tau,K\right) $ for which $K$ is abelian and $%
\omega\left( \tau_{g}(A)^{\ast}\tau_{g}(B)\right) =\omega(A^{\ast}B)$ holds
for all $A,B\in\mathfrak{A}$ and $g\in K$. Let $\varphi_{j}:K\rightarrow K$
be a Borel measurable function such that $\varphi_{j}(gh)=\varphi_{j}(g)%
\varphi_{j}(h)$ for all $g,h\in K$, for $j=1,...,q$. Let $\left\{
\Lambda_{\alpha}\right\} $ be a space-filling net in $K$. For $A,B\in 
\mathfrak{A}$ and $\varepsilon>0$ there then exists an $\alpha_{0}$ such
that 
\begin{align*}
& \left| \frac{1}{\mu(\Lambda_{\alpha_{0}})}\int_{\Lambda_{\alpha_{0}}h}%
\omega\left( A^{\ast}\tau_{\varphi_{1}(g)}(B)\right) ...\omega\left(
A^{\ast}\tau_{\varphi_{q}(g)}(B)\right) dg\right| \\
& >\left| \lim_{\alpha}\frac{1}{\mu(\Lambda_{\alpha})}\int_{\Lambda_{\alpha
}}\omega\left( A^{\ast}\tau_{\varphi_{1}(g)}(B)\right) ...\omega\left(
A^{\ast}\tau_{\varphi_{q}(g)}(B)\right) dg\right| -\varepsilon
\end{align*}
for all $h\in K$. In particular, if $A=B$, then for every $h\in K$%
\begin{equation*}
\left| \omega\left( A^{\ast}\tau_{\varphi_{1}(g)}(A)\right) ...\omega \left(
A^{\ast}\tau_{\varphi_{q}(g)}(A)\right) \right| >\left| \omega\left(
A\right) \right| ^{2q}-\varepsilon 
\end{equation*}
for some $g\in\Lambda_{\alpha_{0}}h$.}

\bigskip

\noindent\textbf{Proof.} We will apply Proposition 4.1 to $\mathfrak{A}_{j}:=%
\mathfrak{A}$, $\omega_{j}:=\omega$ and $\tau_{j,g}:=\tau_{\varphi _{j}(g)}$%
. It is given that $F:K\rightarrow\mathbb{C}:g\mapsto\omega\left(
A\tau_{g}(B)\right) $ is Borel measurable, hence $F\circ\varphi
_{j}:K\rightarrow\mathbb{C}:g\mapsto\omega\left( A\tau_{j,g}(B)\right) $ is
also Borel measurable, while $\tau_{j,g}\circ\tau_{j,h}=\tau_{%
\varphi_{j}(g)}\circ\tau_{\varphi_{j}(h)}=\tau_{\varphi_{j}(g)%
\varphi_{j}(h)}=\tau_{\varphi_{j}(gh)}=\tau_{j,gh}$. Clearly $%
\omega_{j}\left( \tau _{j,g}(A)^{\ast}\tau_{j,g}(B)\right)
=\omega_{j}(A^{\ast}B)$ for all $A,B\in\mathfrak{A}_{j}$ and $\tau_{j,g}(1)=1
$, hence $\left( \mathfrak{A}_{j},\omega_{j},\tau_{j},K\right) $ is a $\ast$%
-dynamical system with the properties required in Proposition 4.1. With $%
A_{j}:=A$ and $B_{j}:=B$, the result now follows from Proposition 4.1. $%
\square$

\bigskip

\noindent\textbf{Corollary 4.3. }\textit{If $\omega(A)>0$, then for every $%
h\in K$ there is a $g\in\Lambda_{\alpha_{0}}h$ such that%
\begin{equation*}
\left| \omega\left( A^{\ast}\tau_{\varphi_{j}(g)}(A)\right) \right| >0 
\end{equation*}
for $j=1,...,q$. (Just take }$\varepsilon<\left| \omega(A)\right| ^{2q}$%
\textit{\ in Theorem 4.2.) $\square$}

\bigskip

For example, since $K$ is abelian, we can take $\varphi_{j}(g)=g^{n_{j}}$
where $n_{j}\in\mathbb{N}$. If $G$ is abelian (or if we just use abelian
notation for $K$), then this says $\varphi_{j}(g)=n_{j}g$, and Corollary 4.3
reduces to the form 
\begin{equation*}
\left| \omega\left( A^{\ast}\tau_{n_{j}g}(A)\right) \right| >0 
\end{equation*}
for $j=1,...,q$.

\end{document}